\documentclass{article}
 \usepackage{amssymb,amsmath,eucal}%,showkeys}

\begin{document}

\input epsf 

\def\boxe{\vphantom{\Bigl |}\boxed}

\def\R {{\mathbb R }}
 \def\C {{\mathbb C }}
  \def\Z{{\mathbb Z}} 
  \def\H{{\mathbb H}}
  \def\N{\mathbb N}
\def\vr{\mathbf {wr}}

\def\frM{\mathfrak M}
\def\SS{\smallskip}

\newcommand{\arcsh}{\mathop{\rm arcsh}\nolimits}
 \newcommand{\sh}{\mathop{\rm sh}\nolimits} 
\newcommand{\ch}{\mathop{\rm ch}\nolimits}
\newcommand{\tha}{\mathop{\rm tanh}\nolimits}

\newcommand{\const}{\mathop{\rm const}\nolimits}

  \newcommand{\Dom}{\mathop{\rm Dom}\nolimits}

% \def\{\addtocounter{punct}{1}{\arabic{sec}.\arabic{punct}.  }}

%\newcounter{punct} [sec]
%\def\punct{\addtocounter{punct}{1}{ \arabic{sec}.\arabic{punct}}}

\def\SL{{\rm SL}}

\def\U{{\rm U}}
\def\O{{\rm O}} 
\def\Sp{{\rm Sp}} 
\def\SO{{\rm SO}}

\def\wh{\widehat}

\def\ov{\overline} 
\def\phi{\varphi} 
\def\epsilon{\varepsilon}
\def\kappa{\varkappa}%don't change to $\kappa$!!!!!!

\def\le{\leqslant} 
\def\ge{\geqslant}

\def\wt{\widetilde}

\def\F{\,\,{}_2F_1} 
\def\FF{\,\,{}_3F_2}

\renewcommand{\Re}{\mathop{\rm Re}\nolimits} 
\renewcommand{\Im}{\mathop{\rm Im}\nolimits}

\newcounter{sec}
 \renewcommand{\theequation}{\arabic{sec}.\arabic{equation}}
\newcounter{punct}[sec]

\newcounter{fact} \def\fact{\addtocounter{fact}{1}{\scc \arabic{fact}}}

\renewcommand{\thepunct}{\thesec.\arabic{punct}}

%\newcounter{punct}
\def\punct{\refstepcounter{punct}{\arabic{sec}.\arabic{punct}.  }}

\def\cL{\mathcal L} 
\def\cM{\mathcal M}
\def\cH{\mathcal H}
\def\cK{\mathcal K}
\def\cE{\mathcal E}
\def\cD{\mathcal D}

\begin{center}
{\Large\bf

Perturbations of Jacobi polynomials
and piece-wise hypergeometric orthogonal systems
}

\medskip

\large

\sc Neretin Yu.A.%
\footnote{Supported by grant NWO--047.017.015}

\end{center}

\bigskip

{\small We construct a family 
of  noncomplete
orthogonal systems of functions on the ray $[0,\infty]$;
these systems
depend on 3 real parameters $\alpha$, $\beta$, $\theta$.
Elements of a system are piece-wise hypergeometric functions,
having a singularity at $x=1$. For
$\theta=0$ these functions vanish on $[1,\infty)$
and our system is reduced to the Jacobi polynomials
$P_n^{\alpha,\beta}$ on the segment $[0,1]$.
In a general case, our functions can be considered
as an interpretation of $P_{n+\theta}^{\alpha,\beta}$.
Our functions are solutions of some exotic
Sturm--Liouville boundary problem for
 the hypergeometric differential
operator. We find the spectral measure for this problem. 
}

\bigskip

{\bf \large 1. Formulation of result}

\medskip

\addtocounter{sec}{1}

Results of the paper are formulated 
in Subsections
 1.1--\ref{l-expansion}.
Next, in \ref{discussion-1}--\ref{attempt},
we discuss existing and hypotetical relations
of our phenomenon with some other
mathematical topics.
 
\SS
%%%%%%%%%%%%%%%%%%%%%%%%%%%%%%%%%%%%%%%%%%%%

{\bf \punct Jacobi polynomials. Preliminaries.
\label{l-jacobi}}
Recall that the Jacobi polynomials
$P_n^{\alpha,\beta}$ are the polynomials
on the segment $[-1,1]$ orthogonal
with respect to the inner product
\begin{equation}
\langle f,g\rangle=
\int_{-1}^1 f(y)\ov{g(y)}(1-y)^{\alpha}(1+y)^{\beta}\,dy
,\qquad \alpha>-1,\beta>-1
\label{jacobi-product}
\end{equation}
These polynomials are given by  explicit formulae,
(see \cite{HTF2}, 10.8(16)),
\begin{align}
P_n^{\alpha,\beta}(y)&=
\frac{\Gamma(n+\alpha+1)}{\Gamma(\alpha+1)n!}
\,\,
F\left[
\begin{matrix}
-n,n+\alpha+\beta+1\\
\alpha+1
\end{matrix}
;\frac{1-y}2\right]
\label{jacobi-1}
\\&=
\frac{(-1)^n\Gamma(n+\beta+1)}{\Gamma(\beta+1) n!}
\,\,
F\left[
\begin{matrix}-n,n+\alpha+\beta+1\\
 \beta+1
 \end{matrix};
  \frac{1+y}2\right]
  \label{jacobi-2}
  \\ &=
\frac{(-1)^n\Gamma(n+\beta+1)}{\Gamma(\beta+1)n!}
\,\,
\Bigl(\frac{1-y}2\Bigr)^{-\alpha}
F\left[
\begin{matrix}
n+\beta+1,-\alpha-n
\\ \beta+1
\end{matrix};   \frac{1-y}2\right]
\label{jacobi-3}
\end{align}

Here $F={}_2F_1$ is the Gauss hypergeometric function,
$$
F\bigl[a,b;c;x\bigr]=
F
\left[\begin{matrix} a,b \\ c\end{matrix};x\right]:=
\sum_{j=0}^\infty \frac{(a)_k(b)_k}{(c)_k} x^k
$$
and $(a)_k:=a(a+1)\dots(a+k-1)$ is the Pochhammer symbol.
 
The expressions (\ref{jacobi-1}), (\ref{jacobi-2})
are polynomials since $(-n)_k=0$ for $k>n$.
The last expression (\ref{jacobi-3}) is a series,
it can be obtained from (\ref{jacobi-2}) by the transformation
(see \cite{HTF1}(2.1.22-23)), 
\begin{equation} 
F\bigl[a,b;c;x\bigr]=
(1-x)^{c-a-b}F\bigl[c-a,c-b;c;x\bigr]
\label{bolza}
\end{equation}

Norms of the Jacobi polynomials with respect
to the inner product
(\ref{jacobi-product})
are given by\begin{equation}
\|P_n^{\alpha,\beta}\|^2=
\langle P_n^{\alpha,\beta}, P_n^{\alpha,\beta}\rangle=
\frac{2^{\alpha+\beta+1}\Gamma(n+\alpha+1)\Gamma(n+\beta+1)}
 {(2n+\alpha+\beta+1)\, n!\,\Gamma(n+\alpha+\beta+1)}
 \label{jacobi-norms}
 \end{equation}

 The Jacobi polynomials are the eigen-functions of
 the differential operator
 \begin{equation}
 D:=(1-y^2) \frac {d^2}{dy^2}
 +\bigl[\beta-\alpha-(\alpha+\beta+2)y)\bigr]\frac d{dy}
 \label{jacobi-operator}
 \end{equation}
 Precisely,
 $$
 D P_n^{\alpha,\beta}=-n(n+\alpha+\beta+1) P_n^{\alpha,\beta}
 $$ 
 
%%%%%%%%%%%%%%%%%%%%%%%%%%%%%%%%%%%%%%%%%%%%%%%%

 {\bf\punct Piece-wise hypergeometric
orthogonal systems%
 \label{l-result}.}
 Now, fix $\theta\in\C$ such that 
 \begin{equation}
 0\le \Re\theta<1
 \label{restriction-theta}
 \end{equation}
 Also, fix $\alpha$, $\beta\in\C$ such that
 \begin{equation}
 -1<\Re \alpha < 1,\,\,\alpha\ne 0,\qquad \Re \beta > - 1
 \label{conditions-1}
 \end{equation} 
 
  Consider the space of functions on the half-line
$x>0$ equipped with the bilinear scalar product
 \begin{equation}
 \{ f,g\}=
 \int\limits_0^1 f(x) g(x) (1-x)^\alpha x^\beta dx+
  \frac{\sin(\alpha+\theta)\pi}
 {\sin \theta\pi}
 \int\limits_1^\infty f(x) g(x) (x-1)^\alpha x^{\beta} {dx}
\label{bilinear-product}
\end{equation} 
 
Denote by $H(x)$ the Heaviside function
$$
H(x)=
\left\{
\begin{aligned}
1,\qquad x>0;\\
0, \qquad x<0
\end{aligned}
\right.
$$
Let $p\in \C$ ranges in the set
\begin{equation}
p-\theta\in\Z,\qquad \Re(2p+\alpha+\beta+1)>0 ,
\quad  1+p+\alpha\ne 0
\label{conditions-2}
\end{equation}

Define the piece-wise hypergeometric
  functions $\Phi_p(x)$ on the half-line $[0,\infty)$
by
\begin{multline}
\Phi_p(x)=
\frac{\Gamma(2p+\alpha+\beta+2)}
{\Gamma(\beta+1)}
F\left[
\begin{matrix}
-p,p+\alpha+\beta+1\\
\beta+1
\end{matrix}; x
\right]\, H(1-x)
+\\+
\frac{\Gamma(1+p+\alpha)}{\Gamma(-p)}
F\left[
\begin{matrix}
p+\alpha+1, p+\alpha+\beta+1\\
2p+\alpha+\beta+2
\end{matrix}; \frac 1x
\right]
x^{-\alpha-\beta-p-1}
 H(x-1)
\label{basis}
\end{multline}

{\sc Theorem 1.} {\it The functions $\Phi_p$ are orthogonal
with respect 
to the symmetric bilinear form (\ref{bilinear-product}),\begin{equation}
\{ \Phi_p,\Phi_q \}=0 \qquad \text{for $p\ne q$}
\label{orthogonality}
\end{equation}
and}
\begin{equation}
\{ \Phi_p,\Phi_p \}=
\frac{\Gamma^2(2p+\alpha+\beta+2)\Gamma(1+p+\alpha)\Gamma(p+1)}
{(2p+\alpha+\beta+1)\Gamma(p+\beta+1)\Gamma(p+\alpha+\beta+1)}
\label{norms}
\end{equation}

We can consider also the Hermitian inner product
 \begin{equation}
 \langle f,g\rangle=
 \int\limits_0^1 f(x)\ov{ g(x)} (1-x)^\alpha x^\beta dx+
  \frac{\sin(\alpha+\theta)\pi}
 {\sin \theta\pi}
 \int\limits_1^\infty f(x)\ov{ g(x)} (x-1)^\alpha x^{\beta} {dx}
\label{hermitian-product}
\end{equation} 
here we must assume $\theta$, $\alpha$, $\beta\in\R$,
and 
$$
0\le\theta<1, \quad -1<\alpha<1,\qquad \beta>-1,\quad 2p+\alpha+\beta+1>0
$$
By our theorem, the functions $\Phi_p$ are orthogonal with
respect to the inner product (\ref{hermitian-product}).
If the factor
$\sin(\alpha+\theta)\pi/
 \sin (\alpha\pi)$ is positive, then our inner product
 also is positive definite.

 \smallskip
  
 {\sc Remark.} {\it The system $\Phi_p$ 
 is not a basis in our Hilbert space.} 
 
 \smallskip
 
{\bf \punct Comparison with the Jacobi polynomials.%
\label{comparison}} 
Let us show, that our construction is reduced to 
the Jacobi polynomials in the case $\theta=0$.

First, assume $x=(1+y)/2$ in the formulae (\ref{jacobi-product}),
(\ref{jacobi-2}). 
We observe that the first summand
in (\ref{basis}) is a Jacobi polynomial.
The second summand in (\ref{basis}) is 0 since
it contains the factor $\Gamma(-p)^{-1}$. This factor is 0
if $p=0,1,2,\dots$. Thus, for an integer $p$, 
$$
\Phi_p(x)=
\frac{\Gamma(2p+\alpha+\beta)\Gamma(p+1)}
     {\Gamma(p+\beta+1)}
     P_p^{\alpha,\beta}(2x-1)H(1-x)
     $$

Hence for $\theta=0$ our orthogonality relations
 are the orthogonality
relations for the Jacobi polynomials.

\smallskip

{\bf \punct Singular boundary problem.
\label{sbp}}
Consider the differential
 operator
\begin{equation}
D:=x(x-1)\frac {d^2}{dx^2} +
(\beta+1-(\alpha+\beta+2)x)\frac d{dx}
\label{D}
\end{equation}
The  functions $\Phi_p$ satisfy  the equation
\begin{equation}
D\Phi_p=-p(p+\alpha+\beta+1)\Phi_p
\label{eigenfunctions}
\end{equation}
More precisely, the function $\Phi_p$ is given by 
different  Kummer 
solutions of the hypergeometric equation
on the  intervals $(0,1)$, $(1,\infty)$,
see formulae \cite{HTF1}, 2.9(1), (13)).

Let $\alpha\ne 0$.%
\footnote{Otherwise, below we have  logarithmic asymptotics
at $x=1$.}
Now we define a space $\cE$ of functions on $[0,\infty)$.
Its elements are  functions $f(x)$ that are smooth
outside the singular points $x=0$, $x=1$, $x=\infty$;
at the singular points they satisfy the following
boundary conditions 
(a strange element of the problem
is the  condition b); its self-adjointness
is not obvious).

\SS

a) {\it The condition at 0.} A function $f$ is smooth at 0.

\SS

b) {\it The condition at 1.} There are functions 
$u(x)$, $v(x)$ smooth at 1 such that
\begin{equation}
f(x)=
\left\{
\begin{aligned}
&u(x)+v(x)(1-x)^{-\alpha}, \qquad &x<1;\\
&\frac{\sin\theta\pi}{\sin(\alpha+\theta)\pi} 
u(x)+v(x)(x-1)^{-\alpha}
&\qquad x>1
\end{aligned}
\right.
\label{uzhas}
\end{equation}

\SS

c) {\it The condition at $\infty$.}
There is a  function $w(y)$ smooth at zero,
such that
$$f(x)=w(1/x)x^{-\alpha-\beta-r-1}\qquad \text{for large $x$}$$
where $r$ is the minimal possible value of $p$.

\smallskip

{\sc Theorem 2.} 
 a) {\it $\Phi_p\in \cE$.}
 
 b) {\it For $f$, $g\in\cE$,}
 $$
 \{Df,g\}=\{f,Dg\}
 $$
 
 Obviously, this implies the orthogonality relations for $p\ne q$.
 Indeed,
 $$\{D\Phi_p, \Phi_q\} = \{\Phi_p, D \Phi_q\}
 =-p(p+\alpha+\beta+1)\{\Phi_p, \Phi_q\}=
 -q(q+\alpha+\beta+1)\{\Phi_p, \Phi_q\}
 $$
 and hence%
\footnote{Under our conditions for parameters,
$p(p+\alpha+\beta+1)=q(q+\alpha+\beta+1)$ implies $p=q$.} $\{\Phi_p, \Phi_q\}=0$.

 \smallskip

 Denote by $\cH$  the Hilbert space with the inner
product (\ref{hermitian-product}).
Obviously, $\cE\subset \cH$.

\smallskip

{\sc Theorem 3.} {\it The operator $D$ 
is essentially self-adjoint on $\cE$.}

\smallskip

{\sc Remark.} a) We can replace the boundary condition
at $\infty$ by the following: $f(x)=0$ for large $x$.
Thus, our complicated formulation is not necessary.

b) If $\beta\ge 1$, then we can replace the condition
at 0 by the following: $f(x)=0$ at some neighborhood of 0.
 For $\beta<1$ the latter simplifying
variant gives a symmetric, but not-self-adjoint
operator%
\footnote{for a discussion of difference
between symmetry and self-adjointness,
see any text-book on the functional
analysis, for instance, \cite{DS}}.
 Possible self-adjoint conditions are
enumerated by points $\lambda:\mu$ of the real projective line;
they can be given in the form
$$
f(x)=A(\lambda+\mu x^{-\beta})+
x\phi(x)+ x^{-\beta+1}\psi(x)
$$
where
$\phi$, $\psi$ are functions 
 smooth near 0, $A$ ranges in $ \C$.
 The condition given above corresponds to $\mu=0$.
Thus, our requirement of the smoothness
is not an absence of a condition,
it hides a condition for asymptotics.
 
 \smallskip

 {\bf \punct Expansion in eigenfunctions%
 \label{l-expansion}.} Our orthogonal  system 
is not complete; hence our operator has a 
partially continuous
specter. In such a case,
a usual expansion of a function in a series of
the Jacobi polynomials must be replaced
by the eigenfunction expansion of $D$ in spirit
of Weyl and Titshmarsh
(see \cite{DS}).

\smallskip

For $s\in \R$, we define the function $\Psi_s(x)$
 on $[0,\infty)$
given by
\begin{multline}
\Psi_s(x)=
F\left[
\begin{matrix}
\frac{\alpha+\beta+1}2+is,\frac{\alpha+\beta+1}2-is\\
\beta+1
\end{matrix};x\right]H(1-x)
+\\+
\frac{2\Gamma(\beta+1)}{\sin (\theta+\alpha)\pi}
\cdot \Re\Biggl\{
\frac{\Gamma(-2is)\cos(\frac{\alpha+\beta}2+\theta-is)\pi}
{\Gamma(\frac{\alpha+\beta+1}2-is)
\Gamma(\frac{-\alpha+\beta+1}2-is)}
\times\\ \times
F\left[
\begin{matrix}
\frac{\alpha+\beta+1}2+is,\frac{\alpha-\beta+1}2+is\\
1+2is
\end{matrix};\frac 1x\right]
x^{-(\alpha+\beta+1)/2-is}\Biggr\}H(x-1)
\label{eigenfunctions-psi}
\end{multline}

Obviously,
$$\Psi_s(x)=\Psi_{-s}(x)$$

{\sc Remark 1.}
The both summands
of $\Psi_s$ are  solutions of the equation
\begin{equation}
D f=-(\frac14(\alpha+\beta+1)^2+s^2)f
\label{yyy}
\end{equation}
Indeed, the first summand is the same as above,
we substitute $p=-\frac{\alpha+\beta+1}2+is$ to (\ref{basis}).
The hypergeometric function in the second summand is 
a Kummer solution of the equation
(\ref{yyy}). Since the coefficients
of $D$ are real, the complex conjugate function
also is a solution of the same equation.

\smallskip

{\sc Remark 2.} The functions $\Psi_s(x)$satisfy the boundary condition at 
$x=1$.
\smallskip

{\sc Remark 3.} The functions $\Phi_p$ with
$p-\theta\in\Z$ are all the $L^2$-eigenfuctions
for the boundary problem formulated above.
Also, the functions $\Psi_s(x)$ are all the remaining
generalized eigenfunctions of the same boundary problem.
See below Section 4.

\smallskip

This 3 remarks easily imply the explicit 
Plancherel measure for the operator $D$.

Consider the Hilbert space $V$, whose elements are 
pairs $(a(p), F(s))$, where $a(p)$ is a sequence
($p$ ranges in the same set as above), and $F(s)$ is a function
on the half-line $s>0$; the inner product is given by
\begin{multline}
\bigl[ (a,F);(b,G)\bigr]
=
\sum_{p}\frac {a(p)\ov{ b(p)}}{\langle\Phi_p,\Phi_p\rangle}+
\\+
\frac {\sin\theta\pi \sin(\theta+\alpha)\pi }{4\pi\Gamma(\beta+1)^2}
\int_0^\infty
\left|
\frac{\Gamma(\frac{\alpha+\beta+1}2-is)
\Gamma(\frac{-\alpha+\beta+1}2-is)}
{\Gamma(2is)\cos(\frac{\alpha+\beta}2 +\theta-is)\pi}\right|^2
F(s)\ov{G(s)} ds
\end{multline}
               
We define the operator 
$Uf\mapsto (a_p,F(s))$ from
$\cH$ to $V$ by
\begin{align*}
a(p)=\langle f,\Phi_p\rangle_{\cH},\\
F(s)=\langle f,\Psi_s\rangle_{\cH}
\end{align*}
 
{\sc Theorem 4.} {\it The operator $U:\cH\to V$
 is a unitary invertible operator.}
 
 \smallskip
 
 In particular, this theorem implies the
 {\it inversion formula},
\begin{equation} 
U^{-1}=U^*
\label{inversion}
\end{equation}

{\bf \punct Discussion:
shift of the index $n$ for classical orthogonal bases%
\label{discussion-1}.}
Thus, for Jacobi
polynomials
$P_n^{\alpha,\beta}$ there is a
perturbated nonpolynomial orthogonal
system obtained by a shift of the number 
  $n$ by a real $\theta$.
Similar deformations are known
for several other classical orthogonal
systems. In the present time, the
general picture looks as confusing,
and I'll shortly refview some known facts.

\smallskip

a) {\it Meixner system.}
Perturbated systems were discovered
by Vilenkin and Klimyk in    \cite{VK},
see, also,  \cite{VK1}, 
and more details in \cite{GK}.

b) {\it  Laguerre system},
see detailed discussion \cite{Gro} 

c) {\it Meixner--Pollachek system,} see. 
\cite{Ner-preprint}.

d) {\it  Continuous dual Hahn
polynomials},
 a (noncomplete) construction
was obtained in \cite{Ner-double}.

In all the cases enumerated above,
the perturbated systems are
orthogonal 
{\it bases}, indexed by numbers 
$n+\theta$, where $n$ ranges in $\Z$.

\SS

All the such deformations are obtained in the following way.
Allmost all classical orthogonal hypergeometric
systems%
\footnote{The only possible exception is 
the Wilson polynomials 
(the author does not know, are they
appear in this context).}
arise in a natural way
in a detailed consideration
of highest weight and lowest weight
representations of
 $\SL_2(\R)$.
Repeating the same operations with
principal and complementary series
of representations,
we obtain deformed systems%
\footnote{For the Hahn system,
the group $\SL_2(\R)$
does not provide a sufficient collection
of parameters,
nevertheless this method has an
heuristic meaning}.

\SS

All the basic formulae existing
for orthogonal polynomials
(see lists in \cite{HTF2}, \cite{KS},
\cite{NUS}),
 "survive"  for deformed systems;
this confirms the nontriviality
of the phenomen under a disscussion.
Moreover, representation-theoretical
interpretation allows to write
such formulae quite easily.

\SS

Note that the systems a)--c) can be
partially%
\footnote{
Partially, since such bases form "series" 
imitating
 "series of representations".
Constructions imitatiting "complementary series",
can be hardly observed without representation theory.} 
produced by very simple
operations.

a)
Fix real parameters $h$, $\sigma$, $t$.
Consider two orthonormal bases 
$$
e_k(z)=z^n,\qquad f_n(z)=
\left(\frac{z\ch t+\sh t}{z\sh t+\ch t}\right)^n
(z\sh t+\ch t)^{-h+i\sigma} (z^{-1}\sh t+\ch t)^{-1+h+i\sigma}
$$
in the space
 $L^2$ on the circle  $|z|=1$.
Expanding one basis in another one,
we obtain a matrix, whose rows
are orthogonal in the space 
$l_2(\Z)$. 
These rows form a Meixner-type system.

 b), c)  Consider an orthonormal
basis in $L^2$
on $\R$ consistng of functions 
$$
f_n(x)=(1+ix)^{-n-1+h+i\sigma}(1-ix)^{n-h+i\sigma}
$$
$n$ ranges in $\Z$.
Applying the Forier transform,
we obtain a  Laguerre-like 
piece-wise confluent hypergeometric basis.

Considering the Mellin transform
$f_n(x)\mapsto (F^+_n(t),F^-_n(t))$,
$$
 F^+_n(t)=\int_0^\infty f_n(x) x^{it-1}\,dx
,\qquad
F^-_n(t)= \int_0^\infty f_n(-x) x^{it-1}\,dx
$$
of the same functions, we obtain a 
Meixner--Pollachek-like system.
It appears that this system is
an orthonormal basis in a space
of $\C^2$-valued functions on $\R$.

d) Hahn-type system can be obtained from
the same functions 
$f_n$ by applying of the bilateral
index hypergeometric transform
introduced in  \cite{Ner-double}.%
\footnote{Again, there arises a
system of 
$\C^2$-valued functions.
There is a general interesting problem
on existence of theory of vector-valued and matrix-valued
special functions. Two examples are given by us just now.}

\SS

Now, the author does not know 
nor represetentation-theoretical interpretation,
nor a way for "simple production" 
Jacobi-type systems.  %
\footnote{The specter of 
the operator $D$
is the same as of a tensor
of a pair of representations of
$\SL_2(\R)$ that are contained in
principal and complementary
series. Possible, this coincidence
is not a chance.}

\smallskip

%%%%%%%%%%%%%%%%%%%%%%%%%%%%%%%%%%%%%%%%%%%%%%%%%

{\bf \punct Discussion: multi-contour
boundary broblems.%
\label{discussion-2}} 
Boundary problems for systems of contours
with cross-glueing of asymptotics
are not well known.
We intend to explain
some natural origins of their appearence.

\SS

First, consider a model example,
the equation
$$
\frac{\partial^2 }{\partial x\partial y}
f(x,y)
=\lambda f(x,y)
$$where $\lambda$ is a constant.
Since, this equation is similar to 
$$
\Bigl(\frac{\partial^2}{\partial^2 x} +
\frac{\partial^2}{\partial^2 y}\Bigr)\,f(x,y)=
\lambda f(x,y)$$
let us imitate the radial separation
of the variables
for the Laplace operator.

\begin{figure}
\caption{}

\qquad
a) $
\epsfbox{kontur.1}
$
\quad\qquad
b) $
\epsfbox{kontur.2}
$

a) Orbits of the group of hyperbolic rotations

b) The space of orbits in general position

\qquad
c) $
\epsfbox{kontur.3}
$
\qquad\qquad
d) $
\epsfbox{kontur.4}
$

c) Asymptotics of 4 functions.

d) Four contours for the Bessel operator $D$.
\end{figure}

Let $f$  be a smooth compactly supported function
on the plane.
For $r>0$ assume
\begin{equation}
g(r,\mu)=\int_{-\infty}^\infty f(re^t,r e^{-t}) e^{-it\mu} d\mu
\label{povorot}
\end{equation}
Respectively, 
 $$
f(x,y)=\frac 1{2\pi} \int_{-\infty}^\infty g(\sqrt{xy},\mu)
\left(\frac xy\right)^{i\mu/2}
d\mu
$$
Our equation now transforms to
\begin{equation}
\cD g=\lambda g
\label{dlambda}
\end{equation}
where
$$ \cD=
\frac 14
\Bigl(\frac {\partial^2}{\partial r^2 }
+\frac 1r \frac {\partial}{\partial r }
+\frac{\mu^2}{r^2}\bigr)
$$

For a fixed $\mu$,
we obtain a Sturm--Liouville problem
for the operator
$D$.

Next, for a fixed
 $\mu\ne 0$ the function $g(r,\mu)$
has the following asymptotics 
$$
g(r,\mu)\sim A r^{i\mu} + B r^{-i\mu},\qquad
r\to 0, 
$$
where
$$A= \frac \pi 2 \Gamma(-i\mu/2) f(0,0)
+\int_0^\infty (f(x,0)-
f(0,0)e^{-x^2})x^{-i\mu-1}dx
$$
and  $B$ is given with a similar
expression.

In the integral transform  (\ref{povorot})
only the values of $f$ 
in the quadrant $x\ge 0$, $y\ge 0$ take part.
Hence, consider functions 
$f(-x,y)$, $f(-x,-y)$, $f(x,-y)$
and construct 3 more functions
 (\ref{povorot}).
As a result, we obtain  4 functions,
whose asymptotics are shown on 
Fig.1c.   
But it two cases we obtain
the equation (\ref{dlambda})
transformed by
$\lambda\mapsto-\lambda$.
Nethertheless, the same result can be
obtained by the substitution
$r\mapsto ir$ in the Bessel operator $\cD$.
Now, in all the 4 cases
the operator 
 $\cD$ became being the same,
but the argument $r$ now
is contained in different contours.

As a result, we obtain a boundary problem of the
following type
(see Fig. 1.d): 
we have the Bessel operator
  $\cD$ defined on quadruples of
functions, 
each function is defined on
its own contour, and the asymptotics at
zero satisfy certain condition of cross-gluing. 

\SS

Emphasis, that our actions were completely standard.
Namely, our equation is invariant
with respect to
action of the group  $\Gamma$
of hyperbolic rotations  $(x,y)\mapsto (xe^t,ye^{-t})$
of the plane, see Fig. 1a).
We simply use this invariance for the separation
of the variables. 
But the set of orbits of $\Gamma$ 
 in a general position 
on $\R^2$ is disconnected and consists 
of 4 components  (see Fig. 1b).
This produces 4-contour problem.

 \SS

More interesting example is the problem
on a spectral decomposition of a
 $\SL_2(\R)$-invariant Laplace operator
on the torus 
$|z|=1$, $|u|=1$,
\begin{multline*}
\Delta=-(z-u)^2 \frac{\partial^2}
{\partial z \partial u}
+\frac 1u (\wt\theta u+\wt\tau z) (z-u)
\frac \partial{\partial z}
+\frac 1z (\theta z+\tau u) (u-z)
\frac \partial{\partial u}+
\\+
\frac 1{zu} (\wt\theta u+\wt\tau z)
(\theta z+\tau u)
\end{multline*}

\begin{figure}
\caption{}
\qquad a) $
\epsfbox{kontur.5}
$
\qquad\qquad
b) $
\epsfbox{kontur.6}
$

a) The diagonal of the torus is the
singular line of the operator 
$\Delta$. Orbits of  $K$ are circles
parallel to the diagonal. The space
of orbits is a circle, but one of point
of this circle is distiguished.

  b)  After a separation
of variables, we obtain a hypergeometric
operator on the contour
 $\Im x=1/2$ (it can be replaced by  $|x|=1$). 

\qquad c) $
\epsfbox{kontur.7}
$
\qquad \qquad
$
\epsfbox{kontur.8}
$

c) The group $P$ has one fixed
point on the torus 
(it is the corner point on our figure).
The space of orbits in a general 
position consists of two components.

d) A separation of variables
produces a confluent hypergeometric
differential operator(in some cases, a bessel
operator)
 on the pair of contours 
 $0$  $\infty$.
The diagonal of the torus
corresponds to the point 0,
and the boundary of the square
to the point  $\infty$.

\qquad
e) $
\epsfbox{kontur.9}
$
\qquad\qquad
$
\epsfbox{kontur.10}
$

e) The group $H$ has   4 fixed pointson the torus. The set of orbits
in a general position consits of
6 components.

 f)  A separation of variables 
produces a hypergeometric
differential operator
on 6 contours. Its singular points
correspond to the diagonal
of the torus and the pieces
$\epsfbox{kontur.11}$, $\epsfbox{kontur.12}$
of the Fig. e).

It os more reasonable to think,
that we have 3 contours
$(\infty,0)$, $(0,1)$, $(1,\infty)$
  and $\C^2$-valued functions
on each contour. 
\end{figure}

The operator is
self-adjoint in  $L^2$ on the torus%
\footnote{There are several more cases
of self-adjointness in other
functional spaces} 
in the case 
 $$
\Re(\theta+\tau)=1,\,\, \Im\theta=\Im\tau,\,\,
\Re(\wt\theta+\wt\tau)=1,\,\, \Im\wt\theta=\Im\wt\tau
$$

This operator is an interesting,
complicated and not well-undertand
object. It was concidered
in several works ,
\cite{Puk}, \cite{Mol-tensor}, \cite{Ner-double}, 
\cite{Gro}). 
We intend to discuss some
details that are absent in these works.

The group $\SL_2(\R)$
acts on the circle by M\"obious transformations,
and hence it acts on the torus since it is
a product of circles. Consider
an one-parametric
group
 $\Gamma\subset\SL_2(\R)$
and separate variables using 
$\Gamma$ as in the previous example.
 
There are 3 possibility

a) $\Gamma=K$  is a subgroup
of rotations of the circle.  

b) $\Gamma=P$ is  parabolic,
i.e. it is an one-parametric subgroup
having one fixed point on
the circle.

c) $\Gamma=H$ is hyperbolic,
i.e., $\Gamma$ is an one-parametric subgroup
having two  fixed points on the circle.

\SS

Separations of variables corresponding to
these subgroups is reduced respectively
to Forier expansion, Fourier transform,
and Mellin transform (a circle is
a  real projective line, and then we can
Fourier and Mellin transform).

As a result, we obtain 3 variants shown
on Fig.2 (we represent the torus as
a square). 

\SS

Numerous nonstandard boundary
problems of this kind 
(in particular, multi-dimensional) arise in
a natural way %
\footnote{Miller's treatise  \cite{Mil}
contains lists of various separations
of variables
for several classical
partial  differential equations.
Some of such ways produce multi-contour problems,
in the book this is not mentioned.},
in the non-commutative harmonic analysis
(the topic apprently arises to  \cite{Ten}).

Actually, they are not well-studied.
Several one-dimensional problems for
the Legendre equation were 
examined by Molchanov
 \cite{Mol-hyp}, \cite{Mol-tensor},
\cite{Mol2}, \cite{Mol3}
and Faraut \cite{Far}
for obtaining of the Plancherel formula
on rank 1 pseudo-Riemannian symmetric spaces
(in particular, this class of spaces includes
multi-dimensional hyperboloids).

\SS

\SS

{\sc Remark}.
Examples enumerated above
allows to think that reasonable multi-contour
problems has approximately the following form.
We consider the space of all the functions
as a module other the space of smooth functions.
Domain of definiteness of a PDE system
consists of functions, that are contained in
a fixed (explicetly defined) submodule
of our module.
A possibility of explicit solution
of such problems looks as questionable;
but also there are no reasons to think
that roundabot ways are better.

\smallskip
{\bf \punct Discussion. Degree of rigidity of the problem.
\label{l-rigidity}}
There are two classical variants of expansion
of the hypergeometric differential operator
in eigenfunctions. 
One case gives the expansion in the Jacobi polynomials.
Another one gives the Weyl--Olevsky index hypergeometric
transform 
\begin{multline}
g(s)
=\\
=\frac 1{\sqrt{2\pi}\Gamma(b+c)}
\int_0^\infty f(x)  \F\Bigl[
  \begin{matrix} b+is,\, b-is\\b+c\end{matrix}
  ; -x\Bigr] x^{b+c-1} (1+x)^{b+c-1}\,dx
\label{def}
\end{multline}
Remind (see, for instance, \cite{DS}, Chapter XIII),
that this transformation is an unitary
operator
$$
L^2\Bigl(\R_+,x^{b+c-1}(1+x)^{b-c}dx\Bigr)
\to
L^2\Bigl(\R_+,\Bigl|\frac
     {\Gamma(b+is)\Gamma(c+is)}{\Gamma(2is)}\Bigr|^2
\Bigr)
$$

Unitarity condition implies the inversion
formula
(\ref{inversion}).
This transformation is an object interesting
by itself having numerious applications
in harmonic analysis and theory of
special functions
(see \cite{Koo}, \cite{VK1}, \cite{Ner-index},
\cite{Ner-wilson}).

There is the third variant of 
this transformation recentely
obtained in
 \cite{Ner-double},
it corresponds to the contour
%(r)н a(r)(r)a??aaa?a?a ?(r)нтуру
shown on Fig. 2b).

Our Theorem 
 4 can be interpreted as one more
(but strange) analog of expansion
in Jacobi polynomials.
Evidentely
 (see, for instance our subsection
\ref{discussion-2}), there are
other analogs
having natural origins. 

On another hand, we can 
assign arbitrary multiplicities to
the contours
$(0,1)$, $(1,\infty)$, $(\infty,0)$,
  on $\C$;
after this there arise a wide
(and even too wide) freedom 
to invent boundaru conditions
as in
(\ref{uzhas}).

As a model example,
 consider the same operator $D$ defined
in  $L^2(0,1)$ with respect to the same weight
$x^\beta(1-x)^\alpha$. If $\alpha>1$, $\beta>1$,
then the space $\cD(0,1)$
 of compactly supported smooth functions
on $(0,1)$ is a domain of self-adjointness. 
If 
$$-1<\beta<1, \qquad -1<\alpha<1$$
   then the deficiency
indices of $D$ on $\cD(0,1)$ are $(2,2)$.
In fact, the both solutions of the equation
$Df=\lambda f$ are in $L^2$ for all $\lambda$.

Fix $\mu$, $\nu\in\R$.
Let us write the boundary conditions
\begin{align*}
f(x)&=A\bigl[1-\mu\frac{\Gamma(-\beta)}{\Gamma(\beta)}
x^{-\beta}\bigr] + x\phi_0(x)+x^{-\beta+1}\psi_0(x)
\quad &\text{near $x=0$}\\
f(x)&=B\bigl[1+\nu\frac{\Gamma(\alpha)}{\Gamma(-\alpha)}
(1-x)^{-\alpha}\bigr] + x\phi_1(x)+(1-x)^{-\alpha+1}\psi_1(x)
\quad &\text{near $x=1$}
\end{align*}
where $\phi_0$, $\psi_0$ are smooth near $x=0$
and $\phi_1$, $\psi_1$ are smooth near $x=1$.

Then the specter is discrete and $\lambda=-p(\alpha+\beta+1+p)$
is a point of  the specter iff
\begin{multline*}
\frac 1{\Gamma(\beta+p+1)\Gamma(-p-\alpha)}+
\frac \lambda {\Gamma(p+1)\Gamma(-\alpha-\beta-p)}
=\\=
\frac\mu{\Gamma(-p)\Gamma(p+\alpha+\beta+1)}
+\frac{\mu\nu}{\Gamma(-\beta-p)\Gamma(\alpha+p+1)}
\end{multline*}
The equation seems nice, but apparently it is nonsolvable.

In any case, not all the boundary conditions
have the equal rights.

{\bf \punct Discussion. An attempt of an application.%
\label{attempt}}
Denote%
\footnote{We imitate a simple way
to derive De Branges--Wilson integral,
proposed by Koornwinder,
see, for instance, \cite{Ner-wilson}.}
$$\xi_\mu=(1-x)^\mu H(1-x), \qquad\mu\in\R$$
Equating
$$
\langle\xi_\mu,\xi_\nu\rangle=[U\xi_\mu,U\xi_\nu]
$$
we obtain the following identity
\begin{multline}
\pi^{-3}\sin\theta\pi \sin(\theta+\alpha)\pi 
\times\\ \times
\int_0^\infty
\left|
\frac{\Gamma(\frac{\alpha+\beta+1}2+is)
\Gamma(\frac{-\alpha+\beta+1}2+is)
\Gamma(\frac{\alpha+\beta+1}2 +\theta+is)
\Gamma(\frac{-\alpha-\beta+1}2 -\theta+is)}
{\Gamma(2is)
\Gamma(\mu+\frac{\alpha+\beta+3}2+is)
\Gamma(\nu+\frac{\alpha+\beta+3}2+is) 
}\right|^2
 ds
 +\\+\\+
 \frac 1{\pi^2}
 \sin(\theta-\mu)\pi \sin(\theta-\nu)\pi
 \sum_p
 (2p+\alpha+\beta+1)\times
 \\ \times
 \Gamma\begin{bmatrix}p+\alpha+\beta+1, p+\beta+1, -\nu+p,-\mu+p\\
 p+\alpha+\beta+\mu+2,  p+\alpha+\beta+\nu+2,p+\alpha+1,p+1
 \end{bmatrix}
 =\\=
 \frac{\Gamma(\beta+1)\Gamma(\alpha+\mu+\nu+1)}
 {\Gamma(\alpha+\beta+\mu+\nu+2)\Gamma(\mu+1)\Gamma(\nu+1)
 \Gamma(\alpha+\mu+1)\Gamma(\alpha+\nu+1)}
\label{beta}
\end{multline}
This identity is a kind of a beta-integral,
continuous and discrete beta-integrals are well-known,
see \cite{Ask}. Our integral has a mixed 
continuous-discrete form, here an integral
and a {\it countable}
${}_6F_5$-sum are present\footnote{ Analytic continuation
of beta-integrals with
respect to parameters can produce
a finite collection of additional summands
 at the left-hand side (due residues).
Our integral has another 
type}. 

Under the substitution $\mu=\theta-1$
 to  (\ref{beta}), the ${}_6F_5$-sum vanishes,
and we obtain the following beta-integral 
obtained in 
 \cite{Ner-wilson},
$$
\frac1{2\pi}\int_{-\infty}^\infty
\Bigl|
\frac{\prod_{i=1}^3 \Gamma(a_k+is)}
{\Gamma(2is)\Gamma(b+is)}
\Bigr|^2
ds
=
\frac{\Gamma(b-a_1-a_2-a_3)\prod_{1\le k<l\le 3}
\Gamma(a_k+a_l)}
{\prod_{i=1}^3 \Gamma(b-a_k)}
$$

The substitution $\theta=0$,
kills the integral term and we obtain a known  
summation formula of  ${}_5F_4$-type.
In fact, this is the famous
Dougall ${}_5H_5$-formula
(see, for instance,  \cite{AAR}.)
$$
\sum_{n=-\infty}^\infty
\frac{\alpha+n}
{\prod_{j=1}^4\Gamma(a_j+\alpha+n)
   \Gamma(a_j-\alpha-n)}
=
\frac{\sin 2\pi\alpha}{2\pi}
\frac{\Gamma(a_1+a_2+a_3+a_4-3)}
{\prod_{1\le j<k\le 4} \Gamma(a_j+a_k-1)}
$$
where one parameter is killed by the substitution
$a_1=\alpha$.

It is interesting
that our integral does not mojorize the Dougall formula.
Apparepntly, this means that 
our construction must
contain additional parameter or
parameters.

\smallskip
 
%{\sc Remark.} There are several substitutions that allow
%to control the correctness of this identity.

%a) Substituting $\theta=0$, we obtain a standard
%${}_5F_4$-identity, see \cite{AAR},
%\cite{ 

%b) Substitute $\mu=-1$. The the integral term is reduced
%to \cite{Ner1}, (0.1). The series term is reduced to
%a ${}_3F_2$-identity \cite{PBM3}, 7.4.4.27.

%c) Pass to a limit as $\mu\to+\infty$, $\nu\to+\infty$.
%The integral term is reduced to the de Branges integral
%(see \cite{AAR}, ) and the discrete sum is reduced to 
%\cite{PBM3}, 7.4.4.27.
%Substituting $\theta=0$, we obtain a standard
%${}_5F_4$-identity

{\bf \punct Structure of the paper.}
In Sections 2 and 3, we give two proofs of the orthogonality
relations. In Section 3 we also discuss our boundary
problem. In Section  4 we obtain the spectral
decomposition of $D$.

\smallskip

{\bf Acknowledgments.}
I am grateful to V.F.Molchanov for a discussion of this
subject.

\medskip

{\bf \large 2. Calculation}

\addtocounter{sec}{1}
\setcounter{equation}{0}
\setcounter{punct}{0}

\medskip

We use the notation
$$\Gamma
\begin{bmatrix} a_1,\dots, a_k\\b_1,\dots b_l
\end{bmatrix}:=
\frac{\Gamma(a_1)\dots \Gamma(a_k)}
{\Gamma(b_1)\dots\Gamma( b_l)}
$$

{\bf \punct The Mellin transform.%
\label{l-mellin}}
For a function $f$ defined on the semi-line $x>0$,
 its {\it Mellin transform} is defined by the formula 
\begin{equation}
\frM f(s)=\int_0^\infty f(x) x^s dx/x
\label{mellin}
\end{equation}
In the cases that are considered below,  this integral  
converges in some strip $\sigma<\Re s <\tau$.
The inversion formula is 
$$
f(x)=\frac 1{2\pi i}
\int_{\gamma-i\infty}^{\gamma+i\infty}
\frM f(s) x^{-s}ds
$$
there the integration is given an over arbitrary contour
lying in the strip $\sigma<\Re s <\tau$.

The multiplicative convolution $f*g$ is defined
by
 \begin{equation}
 f*g(x)=\int_0^\infty f(x)g(y/x)\,dx/x
 \label{convolution}
 \end{equation}
The Mellin transform 
maps the convolution to the product of functions,
\begin{equation}
\frM [f*g](s)= \frM(s) \cdot \frM g(s)
\label{convolution-theorem}
\end{equation}  
(if $\frM f(s)$, $\frM g(s)$ are defined in the 
common strip).

\smallskip

{\bf \punct A way of proof of orthogonality.}
We  write two explicit functions
$\cK_1(s)$, $\cK_2(s)$ and evaluate their
inverse Mellin transforms $K_1$, $K_2$.
Next, we write the identity
 $$K_1*K_2(1)=\frM^{-1}[\cK_1 \cK_2](1)$$
 and observe that it  coincides
 with the orthogonality
 relations for $\Phi_p$ and $\Phi_q$.
 
 \smallskip

Below, in 
  2.7 we explain an origin
of the functions    $K_1$, $K_2$.
The calculation on formal level
(without following of convergences,
conditions for convolutions theorem etc.)
is performed in the next Subsection
2.3. In Subsections 2.4--2.6,
we follow omited details.

\smallskip

{\bf \punct Evaluation of the convolution.%
\label{l-convolution}}
We use the following Barnes-type integral
\cite{PBM3}, 8.4.49.1)
\begin{multline}
\frac 1{2\pi i}
\Gamma
\begin{bmatrix} c,1-b\\a
\end{bmatrix}
\int_{-i\infty}^{+i\infty}
\Gamma
\begin{bmatrix}
s,a-s\\s+1-b,c-s
\end{bmatrix}x^{-s} ds=
\\=
F\left[
\begin{matrix} a,b\\c\end{matrix}
;x\right] H(1-x)
+\\+
x^{-a}
\Gamma
\begin{bmatrix}
c,1-b\\c-a,1+a-b
\end{bmatrix}
F\left[
\begin{matrix}a, 1+a-c\\ 1+a-b
\end{matrix}
;\frac 1x
\right]H(x-1)
\label{main-integral}
\end{multline}
where $x>0$.
The integrand has two series of poles
$$
s=0,-1,-2,\dots,\qquad s=a,a+1,\dots
$$
The integration is given over an arbitrary contour
lying in the strip $0<\Re s<\Re a$
(such contour separates two series of poles).
The condition of the convergence is $\Re(c-a-b)>-1$.

\SS

{\sc Remark.}
A reference to the tables of integrals
is not necessary, since the
identity can be easily proved 
using the Barnes residue method,
see, for instance, (\cite{Sla}, \cite{Mar}).

\SS

We consider two functions $\cK_1(s)$, $\cK_2(s)$ given by
\begin{align*}
\cK_1(s)=
\Gamma
\begin{bmatrix}
\beta+1, p+\alpha+1\\ \beta
+p+1
\end{bmatrix}
\cdot
\Gamma
\begin{bmatrix}
s,\beta+p+1-s\\s+p+\alpha+1,\beta+1-s
\end{bmatrix}
\\
\cK_2(s):=
\Gamma
\begin{bmatrix}
2q+\alpha+\beta+2,-\alpha-q\\q+\alpha+\beta+1
\end{bmatrix}
\cdot
\Gamma
\begin{bmatrix}
s+\alpha+q, \beta+1-s\\
s, q+\beta+2-s
\end{bmatrix}
\end{align*}
we assume that $p-\theta$, $q-\theta\in \Z$.
Using formula
(\ref{main-integral}), we evaluate their inverse
Mellin transforms,
\begin{multline}
K_1(x):=F\left[
\begin{matrix}
p+\beta+1, -p-\alpha\\
\beta+1
\end{matrix}; x\right] H(1-x)
+\\+
x^{-\beta-p-1}
\Gamma\left[
\begin{matrix}
\beta+1,p+\alpha+1\\-p,2p+\alpha+\beta+2
\end{matrix}\right]\cdot
F\left[
\begin{matrix}
\beta+p+1,p+1\\
2p+\alpha+\beta+2
\end{matrix}
;\frac 1x\right] H(x-1)
\end{multline}

\begin{multline*}
K_2(x)=x^{q+\alpha}
F\left[
\begin{matrix} q+\alpha+\beta+1, q+\alpha+1\\
2q+\alpha+\beta+2
\end{matrix};x\right] H(1-x)
+\\+
x^{-\beta-1}
\Gamma
\begin{bmatrix}
2q+\alpha+\beta+2,-\alpha-q\\q+1, \beta+1
\end{bmatrix}
F\left[
\begin{matrix}
q+\alpha+\beta+1, -q\\\beta+1
\end{matrix}
;\frac 1x\right] H(x-1)
\end{multline*}
 Now we write the identity
 $$K_1*K_2(1)=\frM^{-1}[\cK_1 \cK_2](1)$$
and multiply its both sides by  
 \begin{equation}
 \Gamma
 \begin{bmatrix}
 2p+\alpha+\beta+2, q+1\\
 \beta+1, -\alpha-q
 \end{bmatrix}
\label{U}
 \end{equation}
 We obtain the following identity
 \begin{multline}
 \Gamma
 \begin{bmatrix}
 2p+\alpha+\beta+2, \boxe{q+1}\\
 \beta+1, \boxe{-\alpha-q}
 \end{bmatrix}
 \Gamma
\begin{bmatrix}
2q+\alpha+\beta+2,\boxe{-\alpha-q}\\\boxe{q+1}, \beta+1
\end{bmatrix}
\times \\ \times
 \int_0^1
 F\left[
\begin{matrix}
p+\beta+1, -p-\alpha\\
\beta+1
\end{matrix}; x\right]
 x^{\beta+1} F\left[
\begin{matrix}
q+\alpha+\beta+1, -q\\\beta+1
\end{matrix}
;x\right] \,dx/x+
\label{A1}
  \end{multline}
   \begin{multline} 
+ 
 \Gamma
 \begin{bmatrix}
 \boxe{2p+\alpha+\beta+2}, q+1\\
 \boxe{\beta+1}, -\alpha-q
 \end{bmatrix} 
\Gamma\left[
\begin{matrix}
\boxe{\beta+1},p+\alpha+1\\-p,\boxe{2p+\alpha+\beta+2}
\end{matrix}\right]
\times \\ \times
\int\limits_1^\infty
x^{-\beta-p-1}
F\left[
\begin{matrix}
\beta+p+1,p+1\\
2p+\alpha+\beta+2
\end{matrix}
;\frac 1x\right] x^{-q-\alpha}
F\left[
\begin{matrix} q+\alpha+\beta+1, q+\alpha+1\\
2q+\alpha+\beta+2
\end{matrix};\frac 1x\right] \,\frac{dx}x
\label{A2}
\end{multline}  
 $$\!\!\!\!\!\!\!\!\!\!\!\!\!\!\!\!\!\!\!\!\!\!\!\!\!\!\!\!\!\!\!\!\!
 \!\!\!\!\!\!\!\!\!\!\!\!\!\!\!\!\!\!\!\!\!\!\!\!\!\!\!\!\!\!\!\!\!\!
 \!\!\!\!\!\!\!\!\!\!\!\!\!\!\!\!\!\!\!\!\!\!\!\!\!\!\!\!\!\!\!\!\!
 \!\!\!\!\!\!\!
 \!\!\!\!\!\!\!\!\!\!\!\!\!\!\!\!\!\!\!\!\!\!\!\!\!\!\!\!\!\!\!\!\!
 \!\!\!\!\!\!\!\!\!\!\!\!\!\!\!\!\!\!\!\!\!\!\!\!\!\!\!\!\!\!\!\!\!
 \!\!\!\!\!\!\!\!\!\!\!\!\!\!\!\!\!\!\!\!\!\!\!=$$
 \begin{multline}
 \frac 1{2\pi i}
 \Gamma
 \begin{bmatrix}
 2p+\alpha+\beta+2, q+1\\
 \boxe{\beta+1}, \boxe{-\alpha-q}
 \end{bmatrix} 
\Gamma
\begin{bmatrix}
\boxe{\beta+1}, p+\alpha+1\\
 \beta
+p+1
\end{bmatrix}
\Gamma
\begin{bmatrix}
2q+\alpha+\beta+2,\boxe{-\alpha-q}\\q+\alpha+\beta+1
\end{bmatrix}
\times \\ \times
\int_{-i\infty}^{+i\infty}
\Gamma \begin{bmatrix}
\boxe{s},\beta+p+1-s\\s+p+\alpha+1,\boxe{\beta+1-s}
\end{bmatrix}
\Gamma
\begin{bmatrix}
s+\alpha+q, \boxe{\beta+1-s}\\
\boxe{s}, q+\beta+2-s
\end{bmatrix}
\label{A3}
\end{multline}
{\sc Remark.}
Sixteen boxed  $\Gamma$-factors
are intensionally are not canceled.
One of decisive element of the calculation
is cancelations in the row
(\ref{A3}). 
But trick with changing the
integer parameters
 $m$, $n$ by the shifted real
parameters $m+\theta$, $n+\theta$ 
from 2.7 guarantee this cancelation.
 \hfill $\square$

We must identify this identity with the following
orthogonality
identity for $\Phi_p$, $\Phi_q$.
 \begin{multline}
 \Gamma
 \begin{bmatrix}
 2p+\alpha+\beta+2\\
 \beta+1
 \end{bmatrix}
 \Gamma
\begin{bmatrix}
2q+\alpha+\beta+2\\ \beta+1
\end{bmatrix}
\times \\ \times
 \int_0^1
 F\left[
\begin{matrix}
p+\alpha+\beta+1, -p\\ \beta+1
\end{matrix}
;x\right] 
  F\left[
\begin{matrix}
q+\alpha+\beta+1, -q\\ \beta+1
\end{matrix}
;x\right]x^{\beta}(1-x)^\alpha \,dx
+
\label{B1}
\end{multline}
\begin{multline}
+
\frac{\sin(\alpha+\theta)\pi}{\sin\theta\pi}
\,\cdot\,
\Gamma
\begin{bmatrix} 
1+p+\alpha,1+q+\alpha\\-p,-q
\end{bmatrix}
\times\\ \times
\int\limits_1^\infty
F\left[
\begin{matrix} p+\alpha+\beta+1, p+\alpha+1\\
2p+\alpha+\beta+2
\end{matrix};\frac 1x\right]
F\left[
\begin{matrix} q+\alpha+\beta+1, q+\alpha+1\\
2q+\alpha+\beta+2
\end{matrix};\frac 1x\right]
\times\\ \times
x^{-2\alpha-\beta-p-q-2}(x-1)^\alpha \,{dx}=
\label{B2}
\end{multline}
\begin{equation}
=\frac{\delta_{p-q,0}}
{2p+\alpha+\beta+1}
\Gamma
\begin{bmatrix}
2p+\alpha+\beta+2,2p+\alpha+\beta+2,1+p+\alpha,p+1\\
p+\beta+1,p+\alpha+\beta+1
\end{bmatrix}
\label{B3}
\end{equation}
where $\delta_{p-q,0}$
is the Kronecker symbol.

We identify (\ref{A1})--(\ref{A3}) with (\ref{B1})--(\ref{B3})
line-by-line

\smallskip

1. {\it The summand (\ref{A1}) equals the summand (\ref{B1})}.
We transform the first $F$-factor of the integrand
by   (\ref{bolza}).
 
\smallskip

2. {\it The summand (\ref{A2}) equals the summand (\ref{B2})}.
First, we transform the first $F$-factor of the integrand
by   (\ref{bolza}).Secondly, we apply the reflection formula
$\Gamma(z)\Gamma(1-z)=\pi/\sin(\pi z)$
to the gamma-product in  (\ref{A2}).
$$ 
\Gamma
 \begin{bmatrix}
  q+1,p+\alpha+1\\
  -\alpha-q,-p,
\end{bmatrix}=
\frac{\Gamma(1+\alpha+q)\Gamma(1+\alpha+p)\sin(1+\alpha+q)\pi}
{\Gamma(-p)\Gamma(-q)\sin (1+q)\pi}
$$
Next, we use 
$n:=q-\theta\in\Z$,
$$\frac{
\sin(1+\alpha+q)\pi}
{\sin(1+q)\pi}=
\frac{
\sin(1+\alpha+n+\theta)\pi}
{\sin(1+n+\theta)\pi}
=\frac{(-1)^{n+1}
\sin(\alpha+\theta)\pi}
{(-1)^{n+1}\sin(\theta)\pi}
$$

{\sc Remark.}
After this $\Gamma$-factors in  
(\ref{A1}), (\ref{A2}) transforms to the form
$$
\const\cdot u(p)u(q), \qquad \const\cdot v(p)v(q)
$$
with constants that do not depend on
 $p$, $q$.
This is necessary for
a possibility to interprete
the identity 
(\ref{A1})--(\ref{A3})
as an orthogonality relation
of single type functions  $\Phi_p$, $\Phi_q$.
 Certainly, this was achieved due
a multiplication of our identity by
a   $\Gamma$-factor
(\ref{U}).
But, a priory, a possibility
of such multiplications is not
obvious; as far as understand
this is not predictible beforehand. 
\hfill$\square$ 

\SS

3. {\it Right-hand sides, i.e., (\ref{A3}) and (\ref{B3}).}
 We apply the
Barnes-type 
integral
$$
\frac 1{2\pi i}
\int_{-i\infty}^{i\infty}
\Gamma
\begin{bmatrix}
a+s,b-s\\
c+s,d-s
\end{bmatrix}\,ds=
\Gamma
\begin{bmatrix}
a+b, c+d-a-b-1\\
c+d-1, c-a, d-b
\end{bmatrix}
$$
see \cite{PBM3}, 2.2.1.3%
\footnote{This identity is 
a partial case of (\ref{main-integral}).
We substitute $x=1$ to (\ref{main-integral})
and apply the Gauss summation formula for
$F[a,b;c;1]$.} 
and obtain
\begin{multline}
\frac 1{2\pi i}
\int_{-i\infty}^{i\infty}
\Gamma
\begin{bmatrix}
\beta+1+p-s,\alpha+q+s\\
1+p+\alpha+s,q+\beta+2-s
\end{bmatrix}\,ds
=\\=
\Gamma\begin{bmatrix}
\alpha+\beta+p+q+1,\,\, 1\\
\alpha+\beta+p+q+2, q-p+1, p-q+1
\end{bmatrix}
=\\=
\frac 1
{(\alpha+\beta+p+q+1)\Gamma(q-p+1)\Gamma(p-q+1)}
=\\=
\frac 1
{(\alpha+\beta+p+q+1)} 
\frac{\sin(q-p)\pi}{\pi(q-p)}
\label{delta}
\end{multline}
Since $q-p\in \Z$, the latter expression
is zero if 
$p\ne q$.

\smallskip

{\sc Remark.} This place, finishing
the calculation, can look like misterious.
But this almost predictible under the
point of view proposed in 2.7.
Otherwise, how can the Jacobi polynomials
find a possibility to be orthogonal?
\hfill $\square$

\smallskip

{\bf \punct Convergence of the integrals.%
\label{l-convergence}}
%For definiteness, let $\alpha$, $\beta\in\R$.

\smallskip

{\sc Lemma.} {\it Under our conditions
(\ref{conditions-1}), (\ref{conditions-2})
the integral
\begin{equation}
\int_0^1 \Phi_p(x)^2x^\beta(1-x)^\alpha dx
+
\frac{\sin(\alpha+\theta)\pi}{\sin (\theta\pi)}
\int_1^\infty \Phi_p(x)^2 x^\beta (x-1)^\alpha\, dx
\label{integral-for-convergence}
\end{equation}
is absolutely convergent.}

Since
$$|\Phi_p(x)\Phi_q(x)|\le \frac12( |\Phi_p(x)|^2+|\Phi_q(x)|^2)
$$
this lemma implies also the absolute convergence
of
\begin{equation}
\int\limits_0^1 \Phi_p(x)\Phi_q(x)x^\beta(1-x)^\alpha dx
+
\frac{\sin(\alpha+\theta)\pi}{\sin (\theta\pi)}
\int\limits_1^\infty \Phi_p(x)\Phi_q(x) x^\beta (x-1)^\alpha\, dx
\label{integral-for-convergence-2}
\end{equation}

\smallskip

{\sc Proof.}
To follow the asymptotics, we use one of the Kummer
relations (see \cite{HTF1}, 2.10(1)),
\begin{multline}
F\left[
\begin{matrix}
a,b\\c
\end{matrix};
z\right]=
\frac{\Gamma(c)\Gamma(c-a-b)}
     {\Gamma(c-a)\Gamma(c-b)}
F\left[
\begin{matrix}
a,b\\a+b-c+1
\end{matrix};
1-z\right]  
+\\+
\frac{\Gamma(c)\Gamma(a+b-c)}
     {\Gamma(a)\Gamma(b)}
     (1-z)^{c-a-b}
     F\left[
\begin{matrix}
c-a,c-b\\c-a-b+1
\end{matrix};
1-z\right]   
\label{kummer-relation}
\end{multline}

The function $\Phi_p$ is continuous
at $x=0$, and hence the condition
of the convergence of the integral is $\Re\beta>-1$.

The formula (\ref{kummer-relation}) gives
the following asymptotics of $\Phi_p$ as $x\to 1-0$ 
\begin{equation}
C_1+ C_2(1-x)^{-\alpha}
\label{as}
\end{equation}
For $\Re\alpha>0$ we have $\Phi_p\sim (1-x)^{-\alpha}$,
and the condition of convergence 
of (\ref{integral-for-convergence}) is $\Re\alpha<1$.
For $\Re\alpha<0$, the function $\Phi_p$ has a finite limit
at $1+0$, and the condition of convergence
(\ref{integral-for-convergence}) is $\Re\alpha>-1$.

Considering the right limit at 1, we obtain the same
restrictions for $\alpha$.

Obviously,
$$
\Phi_p(x)\sim x^{-\alpha-\beta-p-1},\qquad x\to\infty
$$
Thus the condition of the convergence is
$\Re(2p+\alpha+\beta+1)>0$.

We also must avoid a pole in (\ref{basis}), and this gives
$\alpha+p+1\ne 0$.
\hfill$\square$

\smallskip

Denote $m=p-\theta$, $n=q-\theta$.

\smallskip

{\sc Lemma.}
{\it For fixed $m$, $n$, the integral 
(\ref{integral-for-convergence-2})
depends holomorphicaly on $\alpha$, $\beta$,
$\theta$ in the allowed domain of parameters.
}

\smallskip

{\sc Proof.}
For each given point $(\alpha_0, \beta_0,\theta_0)$,
the convergence of our integral is uniform in a small
neighborhood of $(\alpha_0, \beta_0,\theta_0)$
  (since our asymptotics are uniform).
It remains to refer to the Morera Theorem
(if each integral over closed contour is 0, then the function
is holomorphic).
\hfill $\square$

\smallskip
{\bf \punct Restrictions necessary 
for our calculation.%
\label{l-restrictions}}
First, we used the Mellin transform, and hence
our functions $K_1$, $K_2$ must be locally integrable.
The unique point of discontinuity is $x=1$
We have
\begin{align*}
K_1(x)&\sim A_1+A_2^\pm(1-x)^{-\alpha},\qquad
                   &x\to 1\pm 0;\\
K_2(x)&\sim B_1+B_2^\pm(x-1)^{\alpha},\qquad
                   &x\to 1\pm 0;		   
\end{align*}
This implies $|\Re\alpha|<1$.
		   
Second, we use the convolution theorem for the Mellin
transform.

The Mellin transform (\ref{mellin}) of $K_1$ absolutely converges in
the strip
$$0<\Re s < \beta+p+1$$
The Mellin transform of $K_2$ absolutely converges in the strip
$$-\alpha-q< \Re  <\beta+1$$
We can apply the convolution theorem 
(\ref{convolution-theorem})
if the following conditions are satisfied
\begin{align*}
&\left.
\begin{aligned}
0< \beta+p+1
\\
0< \beta+\alpha+q+1
\end{aligned}
\right\}&
\qquad\text{ --- nonemptyness of strips}\\
&\left.
\begin{aligned}
0< \beta+1
\\
0<p+q+\alpha+ \beta+1
\end{aligned}
\right\}&
\qquad
\begin{matrix}\text{--- nonemptyness of intersection}\\
              \text{of strips}
\end{matrix}	      
	     \end{align*}
 This domain is nonempty, but it is smaller than the domain
 of convergence of (\ref{integral-for-convergence}).
 But the orthogonality identities
 (\ref{orthogonality}), (\ref{norms})
 have holomorphic left-hand sides and right-hand
 sides. Hence they are valid in the whole
 domain of the convergence.
 
 \smallskip

 {\bf \punct Restrictions for $\theta$.}
 These restrictions (\ref{restriction-theta})
 were not used in proof. In fact, $\theta$ is defined
 up to a shift $\theta\mapsto\theta+1$. 
 This shift preserves the orthogonal system
 $\Phi_p$ but changes enumeration
 of the basic elements. 
 
 \smallskip

  By this reason I'll explain how the functions
  $\cK_1$, $\cK_2$ were written. 

  \smallskip
  
 {\bf \punct Comments.  The origin of the calculations.%
 \label{l-comments}}
Now, we intend to explain the origin of
$K_1$, $K_2$.

The orthogonality relations
for the Jacobi polynomials $P_n^{\alpha,\beta}$
are well known but not self-evident.
Let us trying to prove them using
the technique of Barnes integrals, see \cite{Mar}.

Our problem is  an evaluation
of the integral 
\begin{equation}
\int_0^1
\F\left[
\begin{matrix} -n,n+\alpha+\beta\\ \beta+1 
\end{matrix}; x\right]
\F\left[
\begin{matrix} -m,m+\alpha+\beta\\ \beta+1 \end{matrix}; x\right]
x^\beta(1-x)^\alpha dx
\label{jacobi-integral}
\end{equation}
Denote 
\begin{align}
L_1(x)=
(1-x)^\alpha
\F\left[
\begin{matrix} -m,m+\alpha+\beta\\ \beta+1 
\end{matrix}; x\right] H(1-x)
\nonumber
\\
L_2(x)=
x^{-\beta-1}
\F\left[
\begin{matrix} -n,n+\alpha+\beta\\ \beta+1 
\end{matrix}; \frac 1x\right]H(x-1)+r(x) H(x-1)
\label{L2}
\end{align}
where $r(x)$ is an arbitrary function, but
we will choose it later.

Our integral (\ref{jacobi-integral}) is
the convolution  $L_1*L_2(x)$ at the point $x=1$.
We intend to evaluate it using the Mellin
transform.

The Mellin transform of $L_1$
is (see \cite{PBM3}, 8.4.49.1)
$$
\frM L_1(s)=
\Gamma
\begin{bmatrix}
\beta+1, \alpha+m+1\\ \beta+m+1
\end{bmatrix} \cdot
\Gamma
\begin{bmatrix}
s,  \beta+m+1-s\\
\alpha+m+1+s, \beta+1-s
\end{bmatrix}
$$
Then we find a function of the form (\ref{L2})
in the table of inverse Mellin transforms
(see \cite{PBM3}, 8.4.49.1).%
\footnote{
In fact,  tables of integrals are not necessary here,
since we must write a Barnes integral defining
a given hypergeometric function on $[0,1]$,
and it is more-or-less clear how to do this.}

We can assume
$$
\frM L_2(s)=
\Gamma
\begin{bmatrix} 
2n+\alpha+\beta+2,-\alpha-n\\n+\alpha+\beta+1
\end{bmatrix}
\cdot
\Gamma
\begin{bmatrix}
 \alpha+n+s,\beta+1-s\\s, n+\beta+2-s
\end{bmatrix}
$$
and after this the desired calculation can be 
performed.

After this we change $m\to m+\theta$, $n\to n+\theta$.

\SS

{\bf \punct Comments.
 Evaluation of summands in (\ref{B1}), (\ref{B2}).%
\label{l-summands}}
Our orthogonality relations
contain a sum of two integral over different
intervals. It is ineresting
to evaluate each summand separately. 

Let $p$, $q\in\C$.
  Let us evaluate 
$$%\begin{align*}
X:=\int_0^1 \Phi_p(x)\Phi_q(x) x^\beta(1-x)^\alpha dx
,\qquad
Y:=\int_1^\infty \Phi_p(x)\Phi_q(x) x^\beta(x-1)^\alpha dx,\\
$$%\end{align*}
Denote
\begin{align*}
a(p,q)&:=\frac{1}{p+q+\alpha+\beta+1}
\Gamma\bigl[2p+\alpha+\beta+2,2q+\alpha+\beta+2\bigr]
\\
b(p,q)&:=\Gamma\begin{bmatrix}
q+1,p+\alpha+1\\
p+\beta+1,q+\alpha+\beta+1
\end{bmatrix}
\end{align*} 
We write the equation (\ref{A1})--(\ref{A3})
and the same equation with transposed
$p$, $q$
\begin{align*}
X+\frac{\sin(\alpha+q)\pi}{\sin q\pi} Y
=a(p,q)b(p,q) \frac{\sin(q-p)\pi}{(q-p)\pi}
\\
 X+\frac{\sin(\alpha+p)\pi}{\sin p\pi} Y
=a(p,q)b(q,p) \frac{\sin(q-p)\pi}{(q-p)\pi}
\end{align*}

It is a linear system of equations for 
$X$ and $Y$. Its determinant is
$$
\frac{\sin(\alpha+p)\pi}{\sin p\pi}-
\frac{\sin(\alpha+q)\pi}{\sin q\pi}
=
\frac{\sin\alpha\pi \sin(q-p)\pi}
{\sin p\pi \sin q\pi}
$$ 
Hence
\begin{align}
Y&=a(p,q)\frac {\sin p\pi \sin q\pi}
              {\pi(q-p)\sin\alpha\pi}\Bigl[b(q,p)-b(p,q)\Bigr]
	   \label{Y-ravno}   
	      \\
X&=a(p,q)\frac {\sin p\pi \sin q\pi}
              {\pi(q-p)\sin\alpha\pi}
	      \Bigl[\frac{\sin(\alpha+p)\pi}{\sin p\pi}b(p,q)-
	      \frac{\sin(\alpha+q)\pi}{\sin q\pi}b(q,p)\Bigr]=
\nonumber	      \\
&=\frac {\pi a(p,q)}
              {(q-p)\sin\alpha\pi}
\Biggl[\frac 1{\Gamma
     \bigl[p+\beta+1,q+\alpha+\beta+1,-q,-p-\alpha\bigr]}
   -
   \nonumber\\  &\qquad\qquad\qquad\qquad-
\frac 1{\Gamma
     \bigl[q+\beta+1,p+\alpha+\beta+1,-p,-q-\alpha\bigr]}
     \Bigg]
     \label{X-ravno}	      
\end{align}

{\bf \large 3. The boundary problem}

\addtocounter{sec}{1}
\setcounter{equation}{0}
\setcounter{punct}{0}

\medskip

In this section $\alpha\ne 0$.
% But complex values of
%$\alpha$, $\beta$, $\theta$ are admissible.

\SS

{\bf\punct Symmetry of the boundary problem.%
\label{l-symmetry}}
We consider the hypergeometric differential
operator $D$ given by (\ref{D}) and the boundary
problem for $D$ defined in Subsection \ref{sbp}.
We intend to prove the identity
\begin{equation}
\{Df,g\}=\{f,Dg\},\qquad f,g\in\cE
\label{yyyyyy}
\end{equation}

Let 
$$H=a(x)\frac {d^2}{dx^2}+b(x)\frac d{dx}$$
be a differential operator on $[a,b]$ 
formally symmetric with respect
to a weight $\mu(x)$, i.e.,
for smooth $f$, $g$ that vanish near the ends of the interval,
$$\int_a^b Hf(x)\cdot g(x)\,dx=\int_a^b f(x)\cdot H g(x)\,dx$$
Equivalently, $(a\mu)'=b\mu$.
Then for general $f$, $g$, we have
\begin{align}
\int_a^b Hf(x)\cdot g(x)\,dx&-\int_a^b f(x)\cdot H g(x)\,dx=
\nonumber
\\
&=
\Biggl\{\bigl\{f'(x)g(x)-g'(x)f(x)\bigl\}a(x)\mu(x)\Biggr\}\Biggr|_{a}^b
\label{correction}
\end{align}
We apply this identity to the operator $D$
and to the segment $[a,b]=[0,1-\epsilon]$.
Let on some segment $[1-h,1]$ we have
$$
f(x)=u(x)+(1-x)^{-\alpha}v(x),\qquad g(x)=\wt u(x)+(1-x)^{-\alpha}\wt v(x)
$$
with smooth $u(x)$, $v(x)$.
Then the correcting term (\ref{correction}) is
\begin{multline}
\Biggl\{
\det\begin{pmatrix} 
u'(x)+(1-x)^{-\alpha}v'(x)-\alpha(1-x)^{-\alpha-1}v(x)&
     u(x)+(1-x)^{-\alpha} v(x)\\
\wt u'(x)+(1-x)^{-\alpha}\wt v'(x)-\alpha(1-x)^{-\alpha-1}\wt v(x)&     
     \wt u(x)+(1-x)^{-\alpha}\wt v(x)
\end{pmatrix}
\times
\\
\times
x^{\beta+1}(1-x)^{\alpha+1}\Biggr\}\Biggr|_{x=1-\epsilon}
\end{multline}
The last factor gives the power $\epsilon^{\alpha+1}$;
recall that $-1<\alpha<1$. The summands of the determinant
have powers
$$
1,\quad \epsilon^{-\alpha},\quad \epsilon^{-2\alpha}, \quad
\epsilon^{-\alpha-1}\quad, \epsilon^{-2\alpha-1}
$$
But the term with $\epsilon^{-2\alpha-1}$ in the determinant
 is
$$
\det\begin{pmatrix} 
-\alpha v(x)&
      v(x)\\
-\alpha \wt v(x)&     
      \wt v(x)
\end{pmatrix}=0
$$

Hence the leading term of the determinant
has the order
 $\epsilon^{-\alpha-1}$
and  only this term
 gives a contribution to the limit as $\epsilon\to+0$.
 Finally, 
\begin{equation}
 \lim_{\epsilon\to+0}
 \int_{0}^{1-\epsilon} \bigl(D f(x) g(x)-f(x) Dg(x)\bigr)\, 
  x^\beta(x-1)^\alpha\,dx=
u(1)\wt v(1)- v(1)\wt  u(1)
\label{popravka-sleva}
\end{equation}

For $x>1$, we have
$$
f(x)=\frac{\sin\theta\pi}{\sin(\alpha+\theta)\pi}
u(x)+(1-x)^{-\alpha}v(x),\qquad 
g(x)=\frac{\sin\theta\pi}{\sin(\alpha+\theta)\pi}
\wt u(x)+(1-x)^{-\alpha}\wt v(x)
$$
In a similar way,
we obtain
$$
 \lim_{\epsilon\to+0}
 \frac{\sin(\alpha+\theta)\pi}{\sin\theta\pi}
\int_{1+\epsilon}^\infty \bigl(D f(x) g(x)-f(x) Dg(x)\bigr)\, 
x^\beta(x-1)^\alpha\,dx
=-u(1)\wt v(1)+ v(1)\wt  u(1)
$$
This finishes the proof of the identity (\ref{yyyyyy})

\smallskip

{\bf \punct Verification of
the boundary conditions for 
$\Phi_p$.\label{l-phi-p}} Let us show that $\Phi_p(x)$ satisfy 
the boundary conditions at $x=1$.
It is given by a direct calculation,
below we present its details.

We need in expressions for $\Phi_p$
having the form
\begin{equation}
\Phi_p(x)=
\left\{
\begin{aligned}
A_1(p;x)+B_1(p;x)(1-x)^{-\alpha},\qquad x<1
\\
A_2(p,x)+B_2(p,x)(x-1)^{-\alpha}, \qquad x>1
\end{aligned}
\right.
\label{expansion-phi-p}
\end{equation}

We intend to expand $\Phi_p$ in  power series
at $x=1$,
on the semi-segments $(1-\epsilon,1]$, $[1,1+\epsilon)$.

We use the formula (\ref{kummer-relation})
for the left semi-segment and obtain
\begin{multline}
\Phi_p(x)=\Gamma
\begin{bmatrix} 
2p+\alpha+\beta+2\\
\beta+1
\end{bmatrix}
\Biggl\{
\Gamma
\begin{bmatrix}
\beta+1,-\alpha\\
p+\beta+1,-p-\alpha
\end{bmatrix}
F\left[
\begin{matrix}-p,p+\alpha+\beta+1\\ \alpha+1
\end{matrix};\, 1-x\right]
+\\+
\Gamma
\begin{bmatrix}
\beta+1,\alpha\\
-p,p+\alpha+\beta+1
\end{bmatrix}
F\left[
\begin{matrix} p+\beta+1,-p-\alpha\\
               1-\alpha
	       \end{matrix};\,
	       1-x
\right](1-x)^{-\alpha}\Biggr\}
\label{phi-p-left}
\end{multline}	       
for $x<1$
 
Next we use the identity
\begin{multline}
F
\left[
\begin{matrix}
a,b\\c
\end{matrix}; \frac 1x
\right]
=
\Gamma
\begin{bmatrix}
c,c-a-b\\
c-a,c-b
\end{bmatrix}F
\left[
\begin{matrix}
a,a+1-c\\
a+b+1-c
\end{matrix}; 1-x
\right]
x^a
+\\+
\Gamma
\begin{bmatrix}
c,a+b-c\\
a,b
\end{bmatrix}
F
\left[
\begin{matrix}
c-b,1-b\\
c+1-a-b
\end{matrix};
1-x
\right]
x^a(x-1)^{c-a-b} 
\label{kummer-2}
\end{multline}(this formula is a modified variant
of \cite{HTF1}, 1.10(4).
We obtain
\begin{multline}
\Phi_p(x)=
\Gamma
\begin{bmatrix}
p+\alpha+1,-p
\end{bmatrix}
\times
\\
\times
\Biggl\{
\Gamma
\begin{bmatrix}
2p+\alpha+\beta+2,-\alpha\\p+1,p+\beta+1
\end{bmatrix}
F
\left[
\begin{matrix}
p+\alpha+\beta+1,-p\\
\alpha+1
\end{matrix};
1-x
\right]
+\\+
\Gamma
\begin{bmatrix}
2p+\alpha+\beta+2,\alpha\\
p+\alpha+\beta+1,p+\alpha+1
\end{bmatrix}
F
\left[
\begin{matrix}
p+\beta+1,-p-\alpha\\
1-\alpha
\end{matrix};
1-x
\right]
(x-1)^{-\alpha}
\Biggr\}
\label{phi-p-right}
\end{multline}
for $x>1$.

The expressions (\ref{phi-p-left}),
 (\ref{phi-p-right}) are the desired
 expansions (\ref{expansion-phi-p}).
 We observe, that
 $$
 B_1(p,x)=B_2(p,x);
 \qquad A_1(p,x)/A_2(p,x)=\frac{\sin(p+\alpha)\pi}{\sin p\pi}
 $$
 We have
 \begin{equation}
\frac{\sin(p+\alpha)\pi}{\sin p\pi}=\frac{\sin(\theta+\alpha)\pi}{\sin 
\theta\pi}
\label{sin-sin}
\end{equation}
and this implies our boundary condition.

\smallskip

{\sc Remark} (it will be important below in Subsection
\ref{l-adjoint}).
The property (\ref{sin-sin}) is valid iff $p-\theta\in \Z$.
Indeed, the difference between the left-hand side 
and the right-hand
side is 
$$
\frac{\sin\alpha\pi\sin(\theta-p)\pi}
     {\sin p\pi \sin\theta\pi}
$$     
     
\smallskip

{\bf \punct Another proof of the orthogonality 
relations.\label{l-another}}
	For $p\ne q$, the orthogonality follows from
the symmetry condition (\ref{yyyyyy}).

Let us evaluate
$$X:=\int_0^1 \Phi_p(x)\Phi_q(x)x^\beta (1-x)^\alpha dx$$
We preserve the notation (\ref{expansion-phi-p}).
By formula (\ref{popravka-sleva}),
\begin{multline*}
\lim_{\epsilon\to +0}
\Biggl\{\int_0^{1-\epsilon}
D\Phi_p(x) \cdot\Phi_q(x) x^\beta(1-x)^\alpha dx
-
\int_0^{1-\epsilon}
\Phi_p(x) \cdot D\Phi_q(x) x^\beta(1-x)^\alpha dx\Biggr\}
=\\=
A(p,1)B(q,1)-A(q,1)B(p,1)
\end{multline*}
The constants $A(p,1)$ etc. are the $\Gamma$-coefficients
in (\ref{phi-p-left}) and (\ref{phi-p-right});
thus the right-hand side is known.

Since $\Phi_p$ are the eigenfunctions
(see (\ref{eigenfunctions})), the left hand-side is
$$
\bigl[-p(p+\alpha+\beta+1)+
q(q+\alpha+\beta+1)\bigr]\cdot X=
(q-p)(q+p+\alpha+\beta+1)X
$$
After simple cancellations we obtain the 
expression (\ref{X-ravno}).

In the same way, we obtain
the expresion
(\ref{Y-ravno}) for
$\int_1^\infty$.

\medskip

Now we verify our orthogonality
relations via a direct calculation.
But this again is long.

%'???ai мe м(r)ж?м ?a(r)??a?ai a(r)(r)aн(r)шения
%(r)рa(r)?(r)н "iн(r)сти прямым вычис"ением,
%пр вд  (r)пять д"?нн(r)в тым.

{\bf \large 4. The spectral measure}      

\medskip

\addtocounter{sec}{1}
\setcounter{equation}{0}
\setcounter{punct}{0}

Now we intend to evaluate the spectral
measure for the  operator $D$
in the Hilbert space $\cH$ using
Weyl--Titchmarsh machinery, see \cite{DS}.

\smallskip

To avoid logarithmic asymptotics, we
assume $\alpha\ne 0$, $\beta\ne 0$.

{\bf \punct Eigenfunctions of the adjoint
operator.%\label{l-adjoint}
} Now we intend to discuss
the adjoint operator $D^*$ for $D$.

Denote by $\Dom(A)$ the domain of definiteness
of a linear operator $A$.
Recall that $f\in\cH$
is contained in  $\Dom(D^*)$ if there exists a function
$h\in\cH$ such that for each $g\in\Dom(D)$
we have
$$
\langle f, Dg\rangle =
\langle h, g\rangle
$$
In this case, we claim $h=D^*f$.

Since $D$ is symmetric,
we have
$$\Dom(D^*)\supset \Dom(D)=\cE$$

Description of $\Dom(D^*)$ is not an important question,
really it is necessary only description
of eigenfunctions  of $D^*$.

\smallskip

{\sc Lemma.} {\it Let $\Xi$ be an eigenfunction
of $D^*$. Then $\Xi$ satisfies to the boundary conditions
a), b) at $x=0$ and $x=1$ from \ref{sbp}.}

\smallskip\

{\sc Proof.}
 {\it The condition at $0$.}
  Let $D^*\Xi=\lambda\Xi$, represent $\lambda$ as 
\begin{equation}
\lambda=-p(p+\alpha+\beta+1)
\label{lambda-cherez-p}
\end{equation}
There are two solutions of the hypergeometric equation
$Df=\lambda f$ near $x=0$;
if $\beta$ is not a non-negative integer, then they are given
by
\begin{align}
S_1&=\Phi_p(x)=F
\left[
\begin{matrix}-p,p+\alpha+\beta+1\\ \beta+1
\end{matrix};x\right]
\label{nol-1}\\
S_2&=x^{-\beta}
F\left[
\begin{matrix}
-\beta-p,\alpha+p+1\\ 
1-\beta
\end{matrix};x\right]
\label{nol-2}
\end{align}

If $\beta\ge 1$, the second solution
 is not in $\cH$, and the statement is obvious.
 \footnote{For integer $\beta>0$, this also is valid.} 
 
 Let $-1<\beta<1$. 
 Let $f\in \cE$, i.e. $f$ is smooth near $0$.
 Expand our eigenfunction $\Xi$ as
 $$
 \Xi=u(x)+x^{-\beta} v(x), \qquad u(x),\,v(x)\in C^\infty
 $$
 (in fact, $u(x)$ and $v(x)$ are the hypergeometric
 functions defined from (\ref{nol-1}), (\ref{nol-2})
 up to scalar factors.
 If $\Xi$ is in $\Dom(D^*)$, then
 \begin{equation}
 \langle Df,\Xi\rangle-\langle f, D^*\Xi\rangle=0
 \label{xi-in-domain}
 \end{equation}  Repeating the calculation
 of Subsection \ref{l-symmetry},
 we obtain that this difference is
 $$f(0)v(0)$$
 Since $f$ is arbitrary, then $v(0)=0$.
 But 
 $v(x)=\const\cdot F[ -\beta-p,\alpha+p+1; 1-\beta;x]$,
 we have $\const=0$.
 
 \smallskip
 
 {\it The condition at $x=1$.} A proof is similar.
 A priory, we know that
 $$
 \Xi(x)=
 \left\{
 \begin{aligned}
 u_-(x)+v_-(x)(1-x)^{-\alpha}, \qquad x<1\\
 u_+(x)+v_+(x)(x-1)^{-\alpha}, \qquad x<1
 \end{aligned}
 \right.
 $$ 
 In fact $u_\pm$ and $v_\pm$ are the hypergeometric functions
 in the right-hand sides of (\ref{phi-p-left}), (\ref{phi-p-right})
 up to constant factors.

Let $f\in\cE$, i.e.,
$$
f(x)=
\left\{
\begin{aligned}
 a(x)+b(x)(1-x)^{-\alpha},\qquad x<1
 \\
 \frac{\sin \theta\pi} {\sin (\alpha+\theta)\pi}
 a(x)+b(x)(x-1)^\alpha,\qquad x>1
 \end{aligned}
 \right.
 $$
 here $a(x)$, $b(x)$ are smooth near $x=1$,

 If $\Xi\in\Dom(D^*)$, then the condition
 (\ref{xi-in-domain})
 is satisfied. 
 Repeating the considerations of
 Subsection \ref{l-symmetry}, we obtain that 
 (\ref{xi-in-domain}) is equal to 
 $$
 \bigl(
 a(1)v_-(1)-b(1)u_-(1)
 \bigr)-
 \frac{\sin(\alpha+\theta)\pi}
      {\sin \theta\pi}
   \bigl(\frac {\sin \theta\pi}{\sin(\alpha+\theta)\pi}
   a(1)v_+(1)-b(1)u_+(1)\bigr) $$
 It is zero for all $a(1)$, $b(1)$
 and hence
 $$
 v_-(1)=v_+(1),\qquad u_+(1)=
 \frac {\sin \theta\pi}{\sin(\alpha+\theta)\pi}
 u_-(1)
 $$
 But a priory we know $v_\pm$ and $u_\pm$
 up to constant factors, and hence,
 and this implies our statement.
 
 \smallskip
 
 {\bf\punct $L^2$-eigenfunctions of $D^*$.\label{l-adjoint}}
 
 \smallskip
 
 {\sc Lemma.} {\it If $\Xi\in\cH$ is an eigenfunction
 of $D^*$, then $\Xi=\Phi_q$ with $q\in\theta+\Z$.}
 
 \smallskip
 
 {\sc Proof.} Let $\lambda$ be an eigenvalue, let
 $p$ is given by  
 (\ref{lambda-cherez-p}) with $\Re p>-(\alpha+\beta+1)/2$.
 
 Due the boundary condition at $0$,
 we have 
 \begin{equation}
 \Xi=\const F[-p,p+\alpha+\beta+1;\beta+1;x]
 \label{okolo-nulya}
 \end{equation}
  for
 $x<1$. 
 
 Only one solution of the equation $Df=\lambda f$
 is contained in $L^2$ at infinity, it has the form
$$\const\cdot
F[p+\alpha+1,p+\alpha+\beta+1;2p+\alpha+\beta+2;1/x]x^{-\alpha-\beta-p-1}
$$ 
on $[1,\infty]$. Thus, on the both segments
the eigenfunction $\Xi$ coincides 
with  $\Phi_p$ up to scalar factors.
The gluing condition is (\ref{sin-sin}).
By the last remark of Subsection \ref{l-phi-p},
$p-\theta\in\Z$.

If 
\begin{equation}
\Re p=-(\alpha+\beta+1)/2
\label{remaining-specter}
\end{equation} 
then there is no
$L^2$-eigenfunctions at infinity.
\hfill $\square$

\smallskip

{\bf \punct Self-adjointness.}
By the previous lemma, the equations
$D^*f=\pm if$ have no solution in $\cH$.
This implies the essential self-adjointness
of $D$.

\smallskip

{\bf\punct Specter.} The eigenvalues $\lambda=-p(p+\alpha+\beta+1)$
corresponding to the functions $\Phi_p$ form a discrete
specter. The remaining specter corresponds to
the semi-line (\ref{remaining-specter}), i.e., 
$\lambda\ge (\alpha+\beta+1)^2/4$.

Indeed, in all the other cases, we have
precisely one $L^2$ solution $S_0(x)$ of the differential
equation $Df=\lambda f$ near 0, and precisely
one $L^2$-solution $S_\infty(x)$ near infinity.Hence we can write the Green 
kernel
(i.e., the kernel of resolvent) as it is explained
in \cite{DS}. Thus for such $\lambda$
 the resolvent exists.

\smallskip

{\bf \punct Almost $L^2$-eigenfunctions.}
Let
$$p=-(\alpha+\beta+1)/2+is,\qquad s\in\R$$
and $\lambda$ is given by 
(\ref{lambda-cherez-p}).
 
 {\sc Lemma.}
 {\it The function $\Psi_s$
 given by (\ref{eigenfunctions-psi}) is a unique
 almost $L^2$-solution of the equation
 $D\Xi=\lambda \Xi$.}
 
 \smallskip
 
 {\sc Proof.}
 Near $x=0$ such solution must have the form
 (\ref{okolo-nulya}).
 
 We write the following basis $\Lambda(s,x)$, $\Lambda(-s,x)$
 in  the space of solutions of the equation
 $Df=\lambda f$,
 \begin{equation}
 \Lambda(s,x)=
 F
 \left[
 \begin{matrix} 
 \frac{\alpha+\beta+1}2+is, \frac{\alpha-\beta+1}2+is\\
 1+2is
 \end{matrix};\frac 1x
 \right]x^{-(\alpha+\beta+1)/2-is}
 \label{lambda-u-infty}
 \end{equation}

The both solutions are almost $L^2$. Now we must
satisfy the boundary conditions at
$x=1$.
For this, we expand the 3  solutions
(\ref{okolo-nulya}) and $\Lambda(\pm s; x)$
near the point $x=1$.
It remains to write the gluing conditions at $x=1$.
The calculation is long, its reduced to 
usage of the complement formula for $\Gamma$ and elementary
trigonometry. We omit this.
\hfill $\square$

\smallskip

The formula for the spectral measure follows from the explicit asymptotics
of almost $L^2$-solutions at $\infty$;
this is explained in \cite{DS}.

{\sf Math.Physics group, 
 Institute of Theoretical and Experimental Physics, %
\linebreak
B.Cheremushkinskaya, 25, Moscow 117 259, Russia}

%{\sf Moskva 117259, B.  CHeremushkinskaya 25, ITEF, gruppa matfiziki}

{\tt neretin@mccme.ru

\&
Math.Dept,
University of Vienna,
Nordbergstrasse, 15, Vienna 

\end{document}

\end{document}